\documentclass[12pt]{article}
\input epsf.tex

\usepackage{amsmath}
\usepackage{amsthm}
\usepackage{amsfonts}
\usepackage{amssymb}
\usepackage{graphicx}
\usepackage{latexsym}

\usepackage{amsmath,amsthm,amsfonts,amssymb}

\usepackage{epsfig}

\usepackage{amssymb}
\usepackage[all]{xy}
\xyoption{poly}
\usepackage{fancyhdr}
\usepackage{wrapfig}

\theoremstyle{plain}
\newtheorem{thm}{Theorem}[section]
\newtheorem{prop}[thm]{Proposition}
\newtheorem{lem}[thm]{Lemma}

\theoremstyle{definition}
\newtheorem{defn}{Definition}
\theoremstyle{remark}

\advance \topmargin by -\headheight
\advance \topmargin by -\headsep
\evensidemargin \oddsidemargin
\def\F{{\mathbb F}}

\def\N{{\mathbb N}}

\def\chix{{\raise.5ex\hbox{$\chi$}}}

\def\Z{{\mathbb Z}}

\begin{document}
\title{Weak isomorphisms between Bernoulli shifts}
\author{Lewis Bowen\footnote{email:lpbowen@math.hawaii.edu} \\ University of Hawaii}
\begin{abstract}
In this note, we prove that if $G$ is a countable group that contains a nonabelian free subgroup then every pair of nontrivial Bernoulli shifts over $G$ are weakly isomorphic. 
\end{abstract}
\maketitle

\noindent

\section{Introduction}
This paper is motivated by an old and central problem in measurable dynamics: given two dynamical systems, determine whether they are the ``same'', where ``same'' could have one of many different possible meanings. Here we are primarily interested in measure-conjugacy and weak isomorphism. Let us recall some definitions.

A {\em dynamical system} (or {\em system} for short) is a triple $(G,X,\mu)$ where $(X,\mu)$ is a probability space and $G$ is a group acting by measure-preserving transformations on $(X,\mu)$. We will also call this a {\em dynamical system over $G$}, a {\em $G$-system} or an {\em action of $G$}. In this paper, $G$ will always be a discrete countable group. Given two systems $(G,X,\mu)$ and $(G,Y,\nu)$ and a subset $G$-invariant subset $X' \subset X$, a measurable map $\phi:X' \to Y$ is {\em $G$-equivariant} if $\phi(gx)=g\phi(x)$ for all $x\in X'$ and $g\in G$. If, in addition, $\mu(X')=1$ and $\phi_*\mu=\nu$ then $\phi$ is a {\em factor map} from $(G,X,\mu)$ to $(G,Y,\nu)$. If $\phi$ is also a bijection onto its image and $\phi^{-1}:\phi(X') \to X$ is measurable then $\phi$ is a {\em measure-conjugacy}. Two systems are {\em isomorphic}, i.e., {\em measurably conjugate}, if there exists a measure-conjugacy between them. They are {\em weakly isomorphic} if there exist factor maps $\phi:X \to Y$ and $\psi:Y \to X$.


We are concerned with a special class of systems called Bernoulli shifts. To define them, let $(K,\kappa)$ be a standard Borel probability space. For a discrete countable group $G$, let $K^G = \prod_{g \in G} K$ be the set of all functions $x: G \to K$ with the product Borel structure and let $\kappa^G$ be the product measure on $K^G$. $G$ acts on $K^G$ by $(gx)(f)=x(g^{-1}f)$ for $x \in K^G$ and $g,f \in G$. This action is measure-preserving. The system $(G,K^G,\kappa^G)$ is the {\em Bernoulli shift over $G$ with base $(K,\kappa)$}. It is nontrivial if $\kappa$ is not supported on a single point. 

Our main result is:
\begin{thm}\label{thm:main}
If $G$ contains a nonabelian free group then every pair of nontrivial Bernoulli shifts over $G$ are weakly isomorphic.
\end{thm}
It is not known whether this result is true for all countable nonamenable groups. If $G$ is amenable then it is well-known that there exists a pair of Bernoulli shifts over $G$ that are not weakly isomorphic. This follows from the classification of Bernoulli shifts over amenable $G$ up to isomorphism (and weak-isomorphism) that was completed in [OW87].

Let us review the historical context. It was the problem of trying to classify Bernoulli shifts that motivated Kolmogorov to introduce the mean entropy of a dynamical system over $\Z$ [Ko58, Ko59]. That is, Kolmogorov defined for every system $(\Z,X,\mu)$ a number $h(\Z,X,\mu)$ called the {\em mean entropy} of $(\Z,X,\mu)$ that quantifies, in some sense, how ``random'' the system is. To define it, let $\alpha$ be a countable partition of $X$ into at most countably many measurable sets. The {\em entropy} of $\alpha$ is:
$$H(\alpha) = -\sum_{A \in \alpha} \mu(A)\log(\mu(A)).$$
If $\alpha$ and $\beta$ are any two partitions then their {\em join} is defined by $\alpha\vee \beta=\{A \cap B~|~A \in \alpha, B \in \beta\}$. Let $T:X \to X$ be the automorphism implementing the action of $\Z$ on $X$. The {\em mean entropy} of $\alpha$ is
$$h(\alpha):=\lim_{n\to\infty} \big(2n-1\big)^{-1} H\Big( \bigvee_{i=-n}^n T^i\alpha \Big).$$
The partition $\alpha$ is said to be {\em generating} if the smallest $\sigma$-algebra containing $T^i\alpha$ for all $i\in \Z$ is, up to sets of measure zero, equal to the $\sigma$-algebra of all measurable sets in $X$. Kolmogorov proved that the limit defining $h(\alpha)$ exists. Moreover, if $\alpha$ and $\beta$ are generating partitions with $H(\alpha)+H(\beta)<\infty$ then $h(\alpha)=h(\beta)$. He defined the entropy of the system $(\Z,X,\mu)$ to be this common number. 

However, if the system $(\Z,X,\mu)$ does not have a finite-entropy generating partition $\alpha$ then Kolmogorov did not define its entropy. In [Si59] Sinai extended Kolmogorov's definition by defining, for an arbitrary $\Z$-system $(\Z,X,\mu)$, $h(\Z,X,\mu) := \sup_\alpha h(\alpha)$ where the supremum is over all finite-entropy partitions $\alpha$ of $X$. He showed that this definition accords with Kolmogorov's definition when the system possesses a finite-entropy generating partition. This is because entropy is monotone decreasing under factor maps; i.e., if $(\Z,Y,\nu)$ is a factor of $(\Z,X,\mu)$ then $h(\Z,Y,\nu) \le h(\Z,X,\mu)$. Sinai's definition is now standard. Thus $h(\Z,X,\mu)$ is commonly referred to as Kolmogorov-Sinai entropy.


Kolmogorov computed the entropy of a Bernoulli shift over $\Z$ by showing that $h(\Z,K^Z,\kappa^\Z) = H(\kappa)$ where the later is defined as follows. If there exists a finite or countable set $K' \subset K$ such that $\kappa(K')=1$ then 
$$H(\kappa)=-\sum_{k\in K'} \kappa(\{k\})\log\big( \kappa( \{k\} ) \big).$$
Otherwise $H(\kappa)=+\infty$.

Thus two Bernoulli shifts over $\Z$ with different base measure entropies cannot be measurably conjugate. The converse was proven by D. Ornstein in the groundbreaking papers [Or70a, Or70b]. That is, he proved that if two Bernoulli shifts $(\Z, K^\Z, \kappa^\Z), (\Z,L^\Z,\lambda^\Z)$ are such that $H(\kappa)=H(\lambda)$ then they are isomorphic.

 In [Ki75], Kieffer proved a generalization of the Shannon-McMillan theorem to actions of a countable amenable group $G$. In particular, he extended the definition of mean entropy from $\Z$-systems to $G$-systems. This leads to the generalization of Kolmogorov's theorem to amenable groups.

In the landmark paper [OW87], Ornstein and Weiss extended most results of classical entropy theory from $\Z$-systems to $G$-systems where $G$ is any countable amenable group. In particular, they proved that if two Bernoulli shifts $(G, K^G, \kappa^G)$, $(G,L^G,\lambda^G)$ over a countably infinite amenable group $G$ are such that $H(\kappa)=H(\lambda)$ then they are isomorphic. Thus Bernoulli shifts over $G$ are completely classified by base measure entropy.

Now let us say that a group $G$ is {\em Ornstein} if $H(\kappa)=H(\lambda)$ implies $(G, K^G, \kappa^G)$ is isomorphic to $(G,L^G,\lambda^G)$ whenever $(K,\kappa)$ and $(L,\lambda)$ are standard probability spaces. By the above, all countably infinite amenable groups are Ornstein. In [St75] Stepin proved:
\begin{thm}\label{thm:general}
Any countable group that contains an Ornstein subgroup is itself Ornstein.
\end{thm}
We give a proof of this in section \ref{sec:ornstein}. Stepin's theorem and the lemmas preceding its proof are needed to prove theorem \ref{thm:main}. It is unknown whether every countably infinite group is Ornstein.

In [OW87], Ornstein and Weiss asked if entropy theory can be extended to nonamenable groups and presented a curious example which suggests a negative answer. To be precise, let $U_n=\{1,\ldots,n\}$ and let $u_n$ be the uniform probability measure on $U_n$. They showed that if $G$ is the rank 2 free group then $(G,U_2^G,u_2^G)$ factors onto $(G,U_4^G,u_4^G)$. Their example is reproduced in \S \ref{weak isomorphism} (it is a main ingredient in the proof of theorem \ref{thm:main}). This is surprising because, if $G$ is an amenable group then the entropy of $(G,U_n,u_n^G)$ is $\log(n)$. Since entropy is monotone decreasing under factor maps, it follows that $(G,U_2^G,u_2^G)$ cannot factor onto $(G,U_4^G,u_4^G)$. Ornstein and Weiss' factor map does not contradict these facts because the rank 2 free group is not amenable. 

In [Bo08a], the author showed that Kolmogorov's definition of entropy extends to nonabelian free groups, even though Sinai's definition does not. From this extension and Stepin's result it follows that Bernoulli shifts over free groups are completely classified up to measure-conjugacy by base-measure entropy. In [Bo08b], the author extended these results to sofic groups, a class of groups that contains all residually finite groups.

{\bf Acknowledgements}.
It is a pleasure to thank Benjy Weiss for several ideas concerning the proof of theorem \ref{thm:main}. I would also like to thank Dan Rudolph for showing me the theory of coinduced actions and Russ Lyons for pointing out Adam T\'imar's construction.

\section{Stepin's theorem from co-induced actions}\label{sec:ornstein}

In this section, theorem \ref{thm:general} is proven. First, we need to discuss co-induced actions. These actions have been used in [Da06] in an investigation of spectral properties of ergodic actions of discrete groups, in [Ga05] in orbit equivalence theory and in [DGRS08] in constructing non-Bernoulli CPE actions of amenable groups. The definition is related to but different from the well-known Mackey-Zimmer definition of an induced action [Zi78, Zi84]. 

Fix a countable group $G$ and a subgroup $H<G$.
\begin{defn}
A {\em section} for $G/H$ is a map $\sigma: G/H \to G$ such that $\sigma(gH) \in gH$ for all $g\in G$. Fix such a section $\sigma$. For convenience, we assume $\sigma(H)=e$. Let $\alpha:G \times G/H \to H$ be the cocycle 
$$\alpha(g,c)=\sigma(c)^{-1}g\sigma(g^{-1}c)~ \forall g \in G, ~ c\in G/H.$$
$\alpha$ satisfies the cocycle identity:
$$\alpha(g_1g_2,c) = \alpha(g_1,c)\alpha(g_2,g_1^{-1}c) ~\forall g_1,g_2 \in G, ~c \in G/H.$$
Let $(H,W,\omega)$ be a dynamical system over $H$. Define an action of $G$ on the product space $(W^{G/H}, \omega^{G/H})$ by
$$(gx)(c) = \alpha(g,c)x(g^{-1}c) ~ \forall g \in G, ~ x\in W^{G/H}, ~c \in G/H.$$
A routine calculation shows that $(g_1g_2x)(c)=(g_1(g_2x))(c)$ for all $g_1,g_2\in G, x\in W^{G/H}, c \in G/H$. So this action is well-defined.

The system $(G, W^{G/H},\omega^{G/H})$ is said to be {\em co-induced from the action of $H$ on $(W,\omega)$}. Since $H$ preserves $\omega$ and $\alpha(g,c) \in H$ for all $g \in G, c\in G/H$, it follows that $G$ acts on $(W^{G/H},\omega^{G/H})$ by measure-preserving transformations.
\end{defn}

\begin{defn}\label{defn:coinduced factor}
Let $(H,W_1,\omega_1), (H,W_2,\omega_2)$ be dynamical systems over $H$. Suppose $\phi:W_1 \to W_2$ is a factor map. Define $\Phi:W_1^{G/H} \to W_2^{G/H}$ by
$$\Phi(x)(c)=\phi(x(c)),~ \forall x \in W_1^{G/H}, ~ c \in G/H.$$
$\Phi$ is the {\em factor map coinduced by $\phi$}.
\end{defn}

\begin{lem}
$\Phi$ is a factor map from $(G,W_1^{G/H},\omega_1^{G/H})$ to $(H,W_2^{G/H},\omega_2^{G/H})$. If $\phi$ is a measure-conjugacy then so $\Phi$ is also a measure-conjugacy.
\end{lem}
\begin{proof}
 To check $G$-equivariance, note:
$$\Phi(gx)(c)=\phi(gx(c))=\phi(\alpha(g,c)x(g^{-1}c))=\alpha(g,c)\phi(x(g^{-1}c))=\alpha(g,c)\Phi(x)(g^{-1}c)=\big(g\Phi(x)\big)(c).$$
Since $\Phi$ is a product map and $\phi_*\omega_1=\omega_2$, it follows that $\Phi_*\omega_1^{G/H}=\omega_2^{G/H}$. If $\phi$ is invertible then so is $\Phi$ with inverse defined by:
$$\Phi^{-1}(x)(c)=\phi^{-1}(x(c)), ~\forall x \in W_2^{G/H}, ~c \in G/H.$$
\end{proof}

\begin{lem}\label{lem:extra}
Let $(K,\kappa)$ be a probability space. Let $\big(G, (K^H)^{G/H}, (\kappa^H)^{G/H}\big)$ be coinduced from the Bernoulli shift $(H,K^H,\kappa^H)$. Then $\big(G, (K^H)^{G/H}, (\kappa^H)^{G/H}\big)$ is measurably conjugate to the Bernoulli shift $(G,K^G,\kappa^G)$.
\end{lem}

\begin{proof}
Define $J: K^G \to (K^H)^{G/H}$ by 
$$J(x)(c)(h)=x(\sigma(c)h)  , ~ \forall x \in K^G, ~ c \in G/H, ~ h\in H.$$
We claim that $J$ is a measure-conjugacy from $(G,K^G,\kappa^G)$ to $\big(G, (K^H)^{G/H}, (\kappa^H)^{G/H}\big)$.

To check that $J$ is $G$-equivariant, note that for $g \in G, x \in K^G$ and $c\in G/H$, $J(gx)(c)$ is the map $h \mapsto (gx)(\sigma(c)h)=x(g^{-1}\sigma(c)h)$. On the other hand, $(gJ(x))(c) = \alpha(g,c)J(x)(g^{-1}c)$. Since $J(x)(g^{-1}c)$ is the map $h \mapsto x(\sigma(g^{-1}c)h)$, $(gJ(x))(c) $ is the map
$$ h \mapsto x(\sigma(g^{-1}c)\alpha(g,c)^{-1}h) = x( g^{-1}\sigma(c)h).$$
Thus $J(gx)=gJ(x)$. We leave it to the reader to check that $J_{*} \kappa^{G} =( \kappa^H)^{G/H}$ and 
$$J^{-1}(y)(g)=y(gH)(\alpha(g,gH)) ,~\forall y \in (K^H)^{G/H}, ~g \in G.$$
\end{proof}

We can now prove theorem \ref{thm:general}.
\begin{proof}[Proof of theorem \ref{thm:general}]



Let $G$ be a countable group with an Ornstein subgroup $H$. Let $(K_1,\kappa_1), (K_2,\kappa_2)$ be standard probability spaces. Suppose that $H(\kappa_1)=H(\kappa_2)$. We will show that $(G,K_1^G,\kappa_1^G)$ is measurably conjugate to $(G,K_2^G,\kappa_2^G)$. This implies that $G$ is Ornstein.

Because $H$ is Ornstein, there exists a measure-conjugacy $\phi:K_1^H \to K_2^H$ from the $H$-system $(H,K_1^H, \kappa_1^H)$ to $(H,K_2^H, \kappa_2^H)$. By the previous definition, the coinduced map $\Phi:(K_1^H)^{G/H} \to (K_2^H)^{G/H}$ is a measure-conjugacy from $\big(G, (K_1^H)^{G/H}, (\kappa_1^H)^{G/H}\big)$ to $\big(G, (K_2^H)^{G/H}, (\kappa_2^H)^{G/H}\big)$. By the previous lemma $\big(G, (K_i^H)^{G/H}, (\kappa_i^H)^{G/H}\big)$ is measurably conjugate to $(G,K_i^G,\kappa_i^G)$ for $i=1,2$. This proves the theorem.
\end{proof}


\section{Weak Isomorphisms}\label{weak isomorphism}
In this section, we prove theorem \ref{thm:main}. Let $\F=\langle a,b \rangle$ be a rank 2 free group. We begin with Ornstein-Weiss' example from [OW87].
\begin{lem}
For $n\ge 1$, let $U_n=\{1,\ldots,n\}$ and let $u_n$ be the uniform probability measure on $U_n$. Then $(\F,U_2^\F,u_2^\F)$ factors onto $(\F, U_4^\F,u_4^\F)$.
\end{lem}
\begin{proof}
Identify $U_2$ with the group $\Z/2\Z$. Identify $U_4$ with the group $\Z/2\Z \times \Z/2\Z$. Define $\phi:U_2^\F \to U_4^\F$ by
$$\phi(x)(g)=\big( x(g) + x(ga) , x(g)+x(gb)\big).$$
for $x \in U^\F_2$ and $g \in \F$. One way to check that $\phi_*(u_2^\F)=u_4^\F$ is to note that $\phi$ is a group homomorphism where addition in each of $U_2^\F$ and $U_4^\F$ is defined pointwise. Then $u_2^\F$ is Haar measure on $U_2^\F$ and $u_4^\F$ is Haar measure on $U_4^\F$. A simple exercise reveals that $\phi$ is surjective. Thus $\phi_*(u_2^\F)=u_4^\F$.
\end{proof}

Let $U_\infty=(\Z/2\Z)^\N$ and let $u_\infty=u_2^\N$ be the product measure on $U_\infty$. The next lemma is due to Adam T\'imar. It first appeared in [Ba05].
\begin{lem}\label{lem:timar}
There exists a factor map $\Phi$ from $(\F,U_2^\F, u_2^\F)$ to $(\F,U_\infty^\F, u^\F_\infty)$.
\end{lem}
\begin{proof}
For each $n\ge 0$, identify $U_{2^n}$ with $(\Z/2\Z)^n$. Let $\phi:U^\F_2 \to U^\F_4$ be defined as in the previous lemma. Consider the product map:
$$\phi_n: U^\F_{2^n} \times U^\F_2 \to U^\F_{2^n} \times U^\F_4$$
defined by $\phi_n(x,y)=(x,\phi(y))$. We can identify  $U_{2^{n+1}}$ with $U_{2^n}\times U_2$ by $(i_1,\ldots,i_{n+1}) \leftrightarrow \big((i_1,\ldots,i_{n}),i_{n+1}\big).$ By taking the product, we identify $U^\F_{2^{n+1}}$ with $U^\F_{2^{n}} \times  U_2^\F$. Similarly, we identify $U^\F_{2^{n+2}}$ with $U^\F_{2^{n}} \times  U_4^\F$. Therefore, we may regard $\phi_n$ as a map from $U^\F_{2^{n+1}}$ to $U^\F_{2^{n+2}}$. In this aspect, $\phi_n$ is a factor map from $(G,U^\F_{2^{n+1}}, \kappa_{2^{n+1}})$ to $(G,U^\F_{2^{n+2}}, \kappa_{2^{n+2}})$.

Let $\Phi_n: U_2^\F \to U_{2^{n+2}}^\F$ be the composition $\Phi_n := \phi_{n}\circ \cdots \circ \phi_1 \circ \phi$. For $1 \le m < n \le \infty$ define $\pi_m:(\Z/2\Z)^n \to \Z/2\Z$ by $\pi_m(i_1,\ldots,i_n)=i_m.$ If $n_1, n_2 >m$, $x \in U_2^\F$ and $g \in \F$ then 
$$\pi_m\big( \Phi_{n_1}(x)(g) \big ) = \pi_m\big( \Phi_{n_2}(x)(g) \big ).$$
Therefore, we may define $\Phi:U_2^\F \to U_\infty^\F$ by
$$\pi_m\big(\Phi(x)(g)\big) = \pi_m \big(\Phi_{m+1}(x)(g)\big), ~\forall m \in \N, ~x \in U_2^\F, ~g \in \F.$$
$\Phi$ is the required factor map.
\end{proof}

\begin{lem}\label{lem:sinai}
Let $(K_1,\kappa_1)$ and $(K_2,\kappa_2)$ be standard probability spaces. If $H(\kappa_1) \le H(\kappa_2)$ then there exists a factor map $\Phi$ from $(\F,K_2^\F,\kappa_2^\F)$ onto $(\F,K_1^\F,\kappa_1^\F)$.
\end{lem}
\begin{proof}
By Ya. Sinai's factor theorem [Si62] (also in [Gl03, theorem 20.13, page 360]), there exists a factor map $\phi$ from $(\Z,K_2^\Z,\kappa_2^\Z)$  to $(\Z,K_1^\Z, \kappa_1^\Z)$. By definition \ref{defn:coinduced factor} and lemma \ref{lem:extra} the coinduced map $\Phi$ is a factor map from $(\F,K_2^\F,\kappa_2^\F)$ onto $(\F,K_1^\F,\kappa_1^\F)$. Here we have identified $\Z$ with a subgroup of $\F$. 
\end{proof}

 The construction in the following lemma was communicated to the author by Benjy Weiss [We08].
\begin{lem}
Let $(K,\kappa)$ be a standard probability space with $H(\kappa)>0$. Then there exists a factor map from $(\F,K^\F,\kappa^\F)$ onto $(\F,U_\infty^\F,u_\infty^\F)$.
\end{lem}

\begin{proof}
Suppose $H(\kappa) \ge \log(2)$. By lemma \ref{lem:sinai}, $(\F,K^\F,\kappa^\F)$ factors onto $(\F,U_2^\F,u_2^\F)$ which factors onto $(\F,U_\infty^\F,u_\infty^\F)$ by lemma \ref{lem:timar}.

Now suppose that $q=H(\kappa) \in (0,\log(2))$. By Stepin's theorem \ref{thm:general}, $(\F, K^\F,\kappa^\F)$ is measurably-conjugate to $(\F, L^\F,\lambda^\F)$ where $L$ is the three-element space $L=\{0,1,*\}$ and $\lambda(\{0\})=\lambda(\{1\})=p$ where $p >0$ is defined by
\begin{eqnarray}\label{eqn10}
H(\kappa)=-2p\log(p) - (1-2p)\log(1-2p).
\end{eqnarray}

 We identify the subset $\{0,1\} \subset L$ with the group $\Z/2\Z$. Let $M$ be the disjoint union of $\Z/2\Z \times \Z/2\Z$ with $\{*\}$. Define $\phi:L^\F \to M^\F$ as follows. For $x \in L^\F$ and $g \in \F$, if $x(g)=*$, then set $\phi(x)(g):=*$. Otherwise let
$$\phi(x)(g)=\big( x(g) + x(ga^k), x(g) + x(gb^l)\big) \in \Z/2\Z \times \Z/2\Z$$
where $k, l > 0$ are the smallest positive integers such that $x(ga^k) \in \Z/2\Z$ and $x(gb^l) \in \Z/2\Z$. Note $\phi$ is defined on a set of full measure.

We leave it to the reader to check that $\phi$ is a factor map from $(\F,L^\F,\lambda^\F)$ onto $(\F,M^\F, \mu^\F)$ where $\mu$ is defined by $\mu(\{z\})=p/2$ for all $z\in \Z/2\Z\times \Z/2\Z$ and $\mu(\{*\})=\lambda(\{*\})=1-2p$. So $H(\mu)=H(\lambda)+2p\log(2)=H(\kappa)+2p\log(2)$. 

So we have shown that if $H(\kappa) \in (0,\log(2))$ then $(\F,K^\F,\kappa^\F)$ factors onto the Bernoulli shift with base measure entropy $H(\kappa)+2p\log(2)$ where $p$ is defined by equation \ref{eqn10}. Note that $p$ is an increasing function of $H(\kappa)$. So by composing factor maps together, we obtain that $(\F,K^\F,\kappa^\F)$ factors onto the Bernoulli shift with base measure
$$H_n:=H(\kappa)+2p\log(2)+2p_1\log(2)+2p_2\log(2)+\cdots+2p_n\log(2)$$
where $p_1,\ldots,p_n$ satisfy
$$H(\kappa)+2p\log(2)+\cdots+2p_i\log(2) = -2p_{i+1}\log(p_{i+1}) - (1-2p_{i+1})\log(1-2p_{i+1})$$
for $i\in \{1,\ldots,n-1\}$ and $n>0$ is any integer such that $H_{n-1} \le \log(2)$ (where $H_0=H(\kappa)$). Since $p_{i+1} \ge p_i>0$ for all $i$, we have that, for $n$ large enough, $H_n \ge \log(2)$. Thus, $(\F,K^\F,\kappa^\F)$ factors onto a Bernoulli shift with base measure entropy at least $\log(2)$. By the first paragraph, all such shifts factor onto $(\F,U_\infty^\F,u_\infty^\F)$. So, by composing, $(\F,K^\F,\kappa^\F)$ factors onto $(\F,U_\infty^\F,u_\infty^\F)$ too.

\end{proof}

\begin{prop}
If $(K_1,\kappa_1)$ and $(K_2,\kappa_2)$ are any two standard probability spaces with $H(\kappa_1) \ne 0$ then $(\F, K_1^\F,\kappa_1^\F)$ factors onto $(\F,K_2^\F, \kappa_2^\F)$. 
\end{prop}

\begin{proof}
By the previous lemma, $(\F,K_1^\F,\kappa_1^\F)$ factors onto $(\F,U_\infty^\F,u_\infty^\F)$. By lemma \ref{lem:sinai},  $(\F,U_\infty^\F,u_\infty^\F)$ factors onto $(\F,K_2^\F, \kappa_2^\F)$. 
\end{proof}

We can now prove theorem \ref{thm:main}.

\begin{proof}[Proof of theorem \ref{thm:main}]
Let $G$ be a countable group that contains a subgroup $H$ isomorphic to $\F$. Let $(K_1,\kappa_1)$ and $(K_2,\kappa_2)$ be standard probability spaces with $H(\kappa_1) \ne 0$. By the previous proposition, $(H, K_1^H,\kappa_1^H)$ factors onto $(H,K_2^H, \kappa_2^H)$. By definition \ref{defn:coinduced factor} and lemma \ref{lem:extra}, the coinduced map is a factor map from $(G, K_1^G,\kappa_1^G)$ onto $(G,K_2^G, \kappa_2^G)$. Since $(K_1,\kappa_1)$ and $(K_2,\kappa_2)$ are arbitrary, this proves the theorem.
\end{proof}




\end{document}